\documentclass[12pt]{article}
\usepackage{amsthm}
\usepackage{amsmath,amsfonts,amssymb,txfonts,mathrsfs}

\usepackage[all]{xy}
\title{}
\author{Saurav Bhaumik}

\newtheorem{thm}{Theorem}[section]
\newtheorem{lem}[thm]{Lemma}
\newtheorem{proposition}[thm]{Proposition}
\newtheorem{defin}{Definition}
\newtheorem{corol}[thm]{Corollary}

\newenvironment{lemma}{\begin{lem}\;}{\end{lem}}

\theoremstyle{remark}
\newtheorem{remark}[thm]{\bf Remark}

\newcommand{\spec}{{\rm Spec\;}}

\mathchardef\mhyphen="2D

\newcommand{\gm}{\mathbb G_m}

\title{Characteristic Classes for $GO(2n)$ in \'Etale Cohomology}
\date{}
\begin{document}
\maketitle
\begin{abstract}Let $GO(2n)$ be the general orthogonal group (the group of similitudes) over any algebraically closed field of characteristic $\ne 2$. We determine the \'etale cohomology ring with $\mathbb F_2$ coefficients of the algebraic stack $BGO(2n)$. In the topological category, Y. Holla and N. Nitsure determined the singular cohomology ring of the classifying space $BGO(2n)$ of the complex Lie group $GO(2n)$ in terms of explicit generators and relations. We extend their results to the algebraic category. The chief ingredients in this are (i) an extension to \'etale cohomology of an idea of Totaro, originally used in the context of Chow groups, which allows us to approximate the classifying stack by quasi projective schemes; and (ii) construction of a Gysin sequence for the $\gm$-fibration $BO(2n)\to BGO(2n)$ of algebraic stacks.\end{abstract}
Mathematics Subject Classification: 14F20; 14L30 
\section{Introduction}Let $k$ be an algebraically closed field of characteristic $\ne 2$. The {\bf general orthogonal group}, also known as the {\bf group of similitudes}, is the closed subgroup scheme $GO(n)\subset GL_{n,k}$ whose set of $R$-valued points, for any $k$-algebra $R$, is \[GO(n)(R)=\{A\in GL_n(R):\; \exists\, a\in R^\times,\; {}^tAA=aI_n\}.\] This is a reductive group scheme, since it is reduced and its defining representation on $k^n$ is irreducible. In the topological category, Y. Holla and N. Nitsure \cite{nitin1} determined the characteristic classes for the complex Lie group $GO(n)$ with coefficients in $\mathbb F_2$, i.e. they explicitly determined the singular cohomology ring $H^*(BGO(n);\mathbb F_2)$ in terms of generators and relations. The present note proves a similar result in the algebraic category for the \'etale cohomology for the algebraic stack $BGO(n)$.

There are two basic components in our extension of the methods of \cite{nitin1} to the algebraic category. The first of them is a good functorial description of the cohomology of the stack $BGO(n)$. This can be achieved by extending to \'etale cohomology an idea of Totaro \cite{totarochow}, which he originally used for treating Chow groups. The second component is the Gysin sequence for the $\gm$-torsor $BO(n)\to BGO(n)$ in \'etale cohomology, which is actually a by-product of the first component. Once these points are established, the determination of the ring of characteristic classes proceeds exactly as in \cite{nitin1}.

As was remarked in section 1 of \cite{nitin1}, the group scheme $GO(2n+1)$ is a product $\gm\times S\!O(2n+1)$, so $BGO(2n+1)\cong B\gm\times BS\!O(2n+1)$, and \begin{eqnarray*}H^*_{et}(BGO(2n+1);\mathbb F_2)&\cong& H^*_{et}(B\gm;\mathbb F_2)\otimes H_{et}^*(BS\!O(2n+1);\mathbb F_2)\\
 &=& \mathbb F_2[\lambda]\otimes\mathbb F_2[w_2,\ldots,w_{2n+1}]=\mathbb F_2[\lambda,w_2,\ldots,w_{2n+1}].\end{eqnarray*}It is perhaps well known that $H^*_{et}(BS\!O(2n+1);\mathbb F_2)\cong \mathbb F_2[w_2,\ldots,w_{2n+1}]$. Nonetheless, we sketch a proof of this fact at the end of Section \ref{derivation} for the reader's reference. Thus, the case of interest is the cohomology ring of $BGO(2n)$.

\section{Notations and Preliminaries}
Henceforth, all schemes considered in this note will be quasi-projective over an algebraically closed field $k$ of characteristic $\ne 2$, and all morphisms will be over $k$. For such a scheme $X$, we will denote $H^*_{et}(X;\mathbb F_2)$ simply by $H^*(X)$ in what follows. A morphism $X\to Y$ of schemes will be called {\bf $n$-acyclic}, if the induced map $H^i(Y)\to H^i(X)$ is a bijection for $i\le n$. A scheme $X$ will be called {\bf $n$-acyclic} if the structure morphism $X\to\spec k$ is $n$-acyclic. If a morphism is $n$-acyclic for all $n$, it will be called an {\bf acyclic} morphism. A morphism from a scheme to an algebraic stack $X\to \mathfrak X$ over $k$ will be called {\bf $n$-acyclic} if the map induced in cohomology $H^i(\mathfrak X)\to H^i(X)$ is a bijection for $i\le n$.

For any smooth group scheme $G$ over $k$, $BG$ will denote the algebraic stack of all principal bundles on the category of quasi-projective $k$-schemes, locally trivial in \'etale topology. A {\bf characteristic class} for $G$ in \'etale cohomology with coefficient $\mathbb F_2$ is a natural transformation $|BG|\to H_{et}^*(\;;\mathbb F_2)$ of functors from the category of quasi projective schemes over $k$ to the category of sets, where $|BG|$ is defined by taking $|BG|(X)$ to be the set of isomorphism classes in $BG(X)$ for any quasi projective scheme $X$. All such characteristic classes for $G$ form a ring, whose addition and multiplication come from the addition and cup product structure on $H^*$.

\begin{defin}
\label{geomquot}
Let $G$ be a group scheme over $k$, acting on a $k$-scheme $X$. A {\bf bundle quotient} for this action is a morphism $\phi:X\to Y$ of $k$-schemes, such that $X$ together with the given $G$-action is a principal $G$-bundle on $Y$, locally trivial in the \'etale topology. 
\end{defin}

\begin{remark}
\label{bundlequot}
Clearly, whenever such a bundle quotient exists, it is unique up to a unique isomorphism. The existence of bundle quotient is the same as the representability of the quotient stack $[X/G]$ by a scheme. In what follows, bundle quotients will naturally occur in two different ways:\\
a) If a reductive group $G$ acts on a finite type, affine scheme $X$, and $U\subset X$ is an invariant open subscheme such that for each closed point $x\in U$, the orbit $O(x)$ is closed in $X$, and the induced action of $G$ on $U$ is free, then $U$ admits a bundle quotient by $G$. Indeed, since $G$ is reductive, we have a uniform categorical quotient $\phi:X\to X/G$ with $X/G$ affine and $\phi$ universally submersive (Theorem 1.1, \cite{git}). As $U$ is invariant, its image $\phi(U)\subset X/G$ is open. It follows that $\phi:U\to \phi(U)$ is a geometric quotient. Again, since the action on $U$ is free, $\phi:U\to  \phi(U)$ is a principal $G$-bundle (Proposition 0.9 of \cite{git}).\\
b) An affine algebraic group scheme $G$ has a bundle quotient by any of its reduced closed subgroup schemes $H$, which acts on $G$ by right translations. If $G$ is reduced, hence smooth over $k$, then $G/H$, being a reduced homogeneous space, is smooth over $k$ ($k$ being algebraically closed). As $H$ is reduced, hence smooth over $k$, $G\to G/H$ is \'etale locally trivial.

Note that in each of the above cases, the quotient space is a quasi-projective scheme over $k$. 
\end{remark}

\section{Characteristic Classes for $G$ and the Cohomology of the Stack $BG$}
\label{charstack}
\vskip 1em
There is a natural map $\Theta_G$ from the cohomology of the classifying stack $BG$ to the ring of characteristic classes, which takes a cohomology class $\nu\in H^*(BG)$ to the characteristic class whose value on a principal $G$-bundle $P\to X$ is $f^*\nu$, where $f:X\to BG$ is the classifying morphism of $P$. Our treatment of this map is inspired by the article Totaro \cite{totarochow}, which introduced a similar map (actually the inverse of $\Theta_G$) for chow groups, and proved that it was a bijection. 
\begin{proposition}\label{bijection} When $G$ is a reductive group scheme over $k$, the map $\Theta_G$ is a bijection. \end{proposition}
Before we write down the proof of this theorem, let us recall some facts and make a few observations.

Let $Z\hookrightarrow X$ be a closed immersion of two smooth schemes over $k$, where $X$ has relative dimension $n$ over $k$, and all the irreducible components of $Z$ have codimension at least $(s+1)$ in $X$. Then each point $z\in Z$ has an open neighbourhood $z\in U\subset X$, and there is an \'etale morphism $\pi:U\to\mathbb A_k^n$ such that $Z\cap U=\pi^{-1}(\mathbb A_k^{n-t})$ with $t\ge s+1$. Then the theorem of cohomological purity (Theorem 5.1 of \cite{milne}) shows that $\underline{H}^i_Z(X)=0$ for $i<2s+2$. From the spectral sequence $H^i(Z,\underline{H}^j_Z(X))\Rightarrow H^{i+j}_Z(X)$, which shows that $H^i_Z(X)=0$ for $i<2s+2$, and the following long exact sequence \[\cdots\to H^i_Z(X)\to H^i(X)\to H^i(X-Z)\to H^{i+1}_Z(X)\cdots\]we conclude that the open immersion $(X-Z)\hookrightarrow X$ is $2s$-acyclic. However, we need a little stronger statement, where $Z$ need not be smooth:
\begin{lemma}\label{acyclicreduced}Let $Z\hookrightarrow X$ be a closed immersion of finite type schemes over $k$ with $X$ smooth over $k$. If all the irreducible components of $Z$ have codimension $\ge s+1$ in $X$, then the open immersion $(X-Z)\hookrightarrow X$ is $2s$-acyclic.\end{lemma}
\proof The underlying reduced subscheme $Z_{red}$ has a finite filtration by reduced, closed subschemes \[Z_{red}=Z_0\supset Z_1\supset\cdots \supset Z_\ell,\] where for each $i\ge 1$,  $Z_i$ is the singular locus of $Z_{i-1}$ and $Z_\ell$ is non-empty, smooth over $k$. Then ${\rm codim}(Z_i,X)\ge s+i+1$ for each $i$. As $Z_\ell$ is smooth over $k$, $X-Z_\ell\hookrightarrow X$ is $(2s+2l)$-acyclic by the paragraph preceding the Lemma. Note that $Z_{i-1}-Z_i$ is a smooth over $k$, closed subscheme of $X-Z_i$ of codimension $\ge s+i+1$, so that $X-Z_{i-1}\hookrightarrow X-Z_i$ is $(2s+2i)$-acyclic for each $i\ge 1$, by the same argument. Therefore, $X-Z\hookrightarrow X$ is $2s$-acyclic. $\hfill\Box$

\begin{lemma}\label{acyclicbasechange}Let $P\to X$ be a principal $G$-bundle over $k$ with $G$ reduced and $X$ quasi-projective over $k$. Let $f$ be its classifying morphism $X\to BG$. If $P$ is $n$-acyclic over $k$, then the morphism $f$ is also $n$-acyclic.\end{lemma}
\proof We have to show that $R^if_*\mathbb F_{2,X}=0$ where $\mathbb F_{2,X}$ is the constant sheaf on $X$. Let $\pi:\spec k\to BG$ be the classifying morphism of the trivial bundle $G\to \spec k$. As $\pi$ is smooth and surjective, to show that $R^if_*\mathbb F_{2,X}=0$, it is enough to show that $\pi^*R^if_*\mathbb F_{2,X}=0$. Note that we have a Cartesian square 

\centerline{\xymatrix{
P\ar[d]_{f'}\ar[r] & X\ar[d]^f\\
\spec k\ar[r]^-\pi & BG}}

\noindent By smooth base-change, $\pi^*R^if_*\mathbb F_{2,X}=R^if'_*\mathbb F_{2,P}=H^i(P)$, which is zero by hypothesis. \hfill$\Box$

\vskip 1em
\noindent A morphism $X'\to X$ will be called an {\bf affine space bundle} if $X$ admits an \'etale cover $\{X_i\to X\}_{i\in I}$ such that each base change $X'\times_X X_i\to X_i$ is an affine space over $X_i$.

\begin{lemma}\label{jouanolou}Let $X$ be a quasi-projective scheme over $k$, $E$ a principal $G$-bundle on $X$, where $G$ is reductive over $k$, and $s>0$ an integer. Then there is an affine space bundle $\pi:X'\to X$, an affine space $V$ over $k$ with a linear $G$-action on it, and a closed subset $S\subset V$ such that $G$ acts freely on $(V-S)$ and the pullback $\pi^*E$ is isomorphic as a $G$-bundle over $X'$ to the pullback of the principal $G$-bundle $(V-S)\to (V-S)/G$ by a morphism $f:X'\to (V-S)/G$. Moreover, the quotient space $(V-S)/G$ is a quasi-projective scheme over $k$, and the classifying morphism $(V-S)/G\to BG$ for the principal $G$-bundle $(V-S)\to (V-S)/G$ is $2s$-acyclic. \end{lemma}
\proof First see that if $V$ is an affine space with a linear $G$-action, and if $S$ is a closed subset of $V$ such that the orbit of each closed point in $(V-S)$ is closed in $V$, and $G$ acts freely on $(V-S)$, then by Remark \ref{bundlequot} a), we have a bundle quotient $(V-S)\to (V-S)/G$. Therefore, the first statement follows from Lemma 1.6 of \cite{totarochow} (and the proof of it). The quotient space $(V-S)/G$ is quasi-projective by Remark \ref{bundlequot}. By Lemma \ref{acyclicbasechange}, we only have to show that the inclusion $(V-S)\hookrightarrow V$ is $2s$-acyclic, as $V$ is acyclic over $k$. The proof of Lemma 1.6 of \cite{totarochow} shows that $V$ and $S$ can be so chosen that $S$ has codimension $\ge s+1$ in $V$. \hfill$\Box$

\begin{lemma}\label{example}Let $G$ be any reduced algebraic group over $k$. Then, given any positive integer $s$, there exists a linear representation $V$ of $G$, a closed subset $S$ of $V$ such that $G$ acts freely on $(V-S)$, the bundle quotient $(V-S)\to (V-S)/G$ exists, the total space $(V-S)$ is $2s$-acyclic over $k$, and the quotient space is smooth, quasi-projective over $k$. In particular, the classifying morphism of this principal $G$-bundle is $2s$-acyclic.\end{lemma}
\proof This is precisely Remark 1.4 of \cite{totarochow}. To recall the same, let $G$ be regarded as a closed subgroup scheme of $GL_{n,k}$ for some $n$. For any positive integer $s$, choose a positive integer $p$ such that $(p-n)\ge s$. Let $M^o_{n,p}$ be the scheme of $n\times p$ matrices of rank $n$, regarded as an open subscheme of $M_{n,p}=\mathbb A^{n\times p}$. Recall that for the left action of $GL_{n,k}$ on $M^o_{n,p}$, the bundle quotient is the Grassmannian $Gr(k^p,n)$ of $n$-dimensional quotients of $k^p$, and $M^o_{n,p}\to Gr(k^p,n)$ is the tautological principal $GL_n$-bundle on the Grassmannian. Now, the total space of the associated $GL_n/G$-bundle on $Gr(k^p,n)$ is just the bundle quotient of $M^0_{n,p}$ by the action of $G$, which exists as a smooth, quasi-projective scheme over $k$ because it is a quotient of $GL_{p,k}$ by a smooth, closed subgroup scheme (Remark \ref{bundlequot} b)). Also observe that the codimension of the complement of $M^o_{n,p}$ in $M_{n,p}$ is $\ge (n-p+1)\ge (s+1)$. Then by Lemma \ref{acyclicreduced}, we can take $V=M_{n,p}$, and $(V-S)=M^o_{n,p}$. The last sentence follows from Lemma \ref{acyclicbasechange}.\hfill$\Box$
\vskip 1em

\noindent\emph{Proof of Proposition \ref{bijection}.} The injectivity is immediate: if $\Theta_G(\beta)=0$ for a nonzero $\beta\in H^i(BG)$ then by Lemma \ref{example} we find a principal $G$-bundle $P\to X$, whose classifying morphism $X\to BG$ is $i$-acyclic, yielding a contradiction to the definition of $\Theta_G$.

To prove the surjectivity, starting off from a characteristic class $\nu$ of homogeneous degree $s$, say, we choose any pair $(V,S_V)$, where $V$ is a linear representation of $G$, and $S_V$ is a closed subset of $V$ of codimension $\ge s+1$, such that $G$ acts freely on $(V-S_V)$, so that the classifying morphism $(V-S_V)/G\to BG$ is $2s$-acyclic. The value of $\nu$ on this principal $G$-bundle is the image, under the induced map in cohomology, of a unique element, say $\alpha\in H^s(BG)$. We claim that $\nu=\Theta_G(\alpha)$. The claim will follow immediately from Lemma \ref{jouanolou}, once we show that $\alpha$ is independent of the pair $(V,S_V)$. But this in turn is a consequence of Totaro's argument involving ``independence of $V$ and $S$'', which we recall from the proof of Theorem 1.1 of \cite{totarochow}: let $(W,S_W)$ be another pair such as $(V,S_V)$. Then the bundle quotients $((V-S_V)\times W)/G$ and $(V\times(W-S_W))/G$ exist as vector bundles over $(V-S_V)/G$ and $(W-S_W)/G$, respectively. Both the closed subsets $V\times S_W$ and $S_V\times W$ of $V\oplus W$ have codimensions at least $(s+1)$, and outside each of them, $G$ acts freely, so that in each case the bundle quotient exists. Therefore, there is an invariant open subscheme of $V\oplus W$ which contains both of these open complements and which consists of $G$-stable points. Let us call this open subscheme $V\oplus W-S_{V\oplus W}$. Observe that $S_{V\oplus W}$, being contained in each of the two closed subsets mentioned above, has codimension at least $(s+1)$ in $V\oplus W$. Since $H^i((V-S)/G)$ for $i\le 2s$ does not depend on $S$ as long as the codimension of $S$ is $\ge s+1$, we see that all the arrows in the following commutative diagram are isomorphisms for $i\le 2s$, where the right vertical, the diagonal, and the lower horizontal arrows are induced by the classifying morphisms, and the upper horizontal and the left vertical arrows, by obvious inclusions: 

\centerline{\xymatrix{
H^i((V\oplus W-S_{V\oplus W})/G)\ar[r]\ar[d]& H^i((V\times(W-S_W))/G)= H^i((W-S_W)/G)\\
H^i((V-S_V)/G)=H^i(((V-S_V)\times W)/G)&H^i(BG)\ar[l]\ar[u]\ar[lu]
}}

\noindent This completes the proof.
 $\hfill\Box$
\begin{remark}\label{modification}The independence of $V$ and $S$ as above can also be applied in a little different setup: suppose $P\to X$ is a principal $G$-bundle, and there are two $G$-equivariant pairs $(V,S_V)$ and $(W,S_W)$ with maps from $X$ to the bundle quotients $(V-S_V)/G$ and $(W-S_W)/G$, with properties as required in Lemma \ref{jouanolou}. Then the diagonal map $P\to V\oplus W$ factors through an open inclusion $(V\oplus W-S_{V\oplus W})\hookrightarrow V\oplus W$ as above.\end{remark}

\section{The Gysin Sequence}
\vskip 1em
We begin by observing the following:
\begin{remark}\label{bundlegysin}Let $L\to X$ be a line bundle over $k$. A line bundle being locally trivial in the Zariski topology, the zero section has \emph{pure codimension} in $L$. Therefore, whenever there is a principal $\gm$-bundle with the quotient space, hence the total space as well, smooth over $k$, we have the Gysin sequence for the smooth pair $(L,X)$, where $X$ is regarded as the zero section.\end{remark}
\begin{lemma}\label{associated}Let $1\to N\to G\to G'\to 1$ be an exact sequence of reduced, affine algebraic groups over $k$, and let $E\to X$ be a principal $G$-bundle. Then the morphism $E\to X$ admits a factorization $E\to E'\to X$, where $E\to E'$ is the bundle quotient of $E$ by $N$, and $E'\to X$ is the $G'$-bundle associated to $E\to X$.\end{lemma} 
\proof The associated $G'$-bundle $E'\to X$ exists by \'etale descent of affine morphisms. It is seen locally that the morphism $E\to E'$ is the required bundle quotient. $\hfill\Box$

In the above, let $G'=\gm$, the multiplicative group scheme. Let $P\to Q$ be a principal $G$-bundle, where $Q$ is a smooth, quasi-projective scheme over $k$. Then the bundle quotient of $P$ by the action of $N$ exists as the $\gm$-bundle $Q'\to Q$ associated to $P\to Q$, by Lemma \ref{associated}. We can recognise $Q'\to Q$ as the complement of the zero section of the induced line bundle. The Gysin long exact sequence \[\cdots\to H^i(Q)\to H^i(Q')\to H^{i-1}(Q)\to H^{i+1}(Q)\to\cdots\]exists by Remark \ref{bundlegysin}, and is functorial for maps of the quotient space. 

On the other hand, we have the universal principal $G$-bundle $\spec k\to BG$, whose associated $\gm$-torsor is the 1-morphism $BN\to BG$. Assuming that $G$ is reductive, we will construct a long exact sequence for this $\gm$-torsor: \[\cdots\to H^i(BG)\to H^i(BN)\to H^{i-1}(BG)\to H^{i+1}(BG)\to \cdots\]In order to do this, we observe that by Lemma \ref{example}, given any integer $s>0$, there is a principal $G$-bundle $P\to Q$ with $H^i(P)=0$ for $i\le 2s+2$, so that the classifying morphism $Q\to BG$ is $(2s+2)$-acyclic. Similarly, if $Q'\to Q$ is the associated $\gm$-bundle, then $P\to Q'$ is a principal $N$-bundle, whose classifying morphism $Q'\to BN$ is $(2s+2)$-acyclic. We define the Gysin sequence for $BN\to BG$ by requiring each of the squares in the following diagram to commute, where the vertical isomorphisms are induced by the classifying morphisms for $i\le 2s$:

\centerline{\xymatrix{
\cdots\ar[r]& H^i(BG)\ar[r]\ar[d]_\cong & H^i(BN)\ar[r]\ar[d]_\cong & H^{i-1}(BG)\ar[r]\ar[d]_\cong & H^{i+1}(BG)\ar[r]\ar[d]_\cong &\cdots\\
\cdots\ar[r]& H^i(Q)\ar[r]& H^i(Q')\ar[r]& H^{i-1}(Q)\ar[r] & H^{i+1}(Q)\ar[r] & \cdots
}}

That the Gysin sequence does not depend on the particular principal bundle chosen, but any principal $G$-bundle $P\to Q$ of the type mentioned in the statement of Lemma \ref{example} defines the same Gysin sequence for $BN\to BG$, follows from the ``independence of $(V,S)$ argument'', which was used in the proof of Proposition \ref{bijection}. More generally, by the independence of $(V,S)$ argument and Lemma \ref{jouanolou}, the Gysin sequence for $BN\to BG$ is compatible with the Gysin sequence for principal $G$-bundles over schemes (even when the vertical morphisms are not isomorphisms) by the same independence of $(V,S)$ argument. Therefore we have proved the following proposition.
\begin{proposition}\label{gysin}If $G$ is reductive, there is a Gysin long exact sequence for the $\gm$-torsor $BN\to BG$, which is compatible with Gysin sequence for principal $G$-bundles over schemes in the sense that if $P\to Q$ is a principal $G$-bundle with $Q$ a smooth, quasi-projective scheme over $k$, and $Q'\to Q$ the associated $\gm$-bundle, then  all the squares in the following diagram are commutative, where the vertical maps are induced by the classifying morphisms.

\centerline{\xymatrix{
\cdots\ar[r]& H^i(BG)\ar[r]\ar[d]& H^i(BN)\ar[r]\ar[d]& H^{i-1}(BG)\ar[r]\ar[d]& H^{i+1}(BG)\ar[r]\ar[d]&\cdots\\
\cdots\ar[r]& H^i(Q)\ar[r]& H^i(Q')\ar[r]& H^{i-1}(Q)\ar[r] & H^{i+1}(Q)\ar[r] & \cdots
}}\end{proposition}

\section{Cohomology Calculation}
\label{derivation}
In the paper of Holla and Nitsure \cite{nitin1}, the main ingredient of the determination of the singular cohomology ring $H^*_{sing}(BGO(2n);\mathbb F_2)$ was the Gysin sequence. Other conceptual points that were used their argument are the following :\\
(1) The fact that $H^*_{sing}(BO(n);\mathbb F_2)\cong \mathbb F_2[w_1,\ldots,w_n]$, where $w_i$ is the $i$-th Stiefel-Whitney class;\\
(2) the splitting principle for $O(n)$-bundles;\\
(3) the fact that under the map in cohomology induced by the inclusion $BO(n)\subset BGL_n$, $c_i\mapsto w_i^2$, and \\
(4) the K\"unneth formula, i.e. that the natural map $H^*_{sing}(\mathbb C^\times;\mathbb F_2)\otimes H^*_{sing}(X;\mathbb F_2)\to H^*_{sing}(\mathbb C^\times\times X;\mathbb F_2)$ is an isomorphism. 

The Gysin sequence for $BO(2n)\to BG(2n)$ in \'etale cohomology is already established in Proposition \ref{gysin}. Among the other conceptual points, an \'etale cohomology version of (1) is available in \cite{hassewit}. (2) The splitting principle in \'etale cohomology is the `Claim' at page 180 of \cite{helen}, while an \'etale cohomology version of (3) appears in the same paper. We also have the K\"unneth formula in \'etale cohomology for the following special case: The natural map $H^*_{et}(X;\mathbb F_2)\otimes H^*_{et}(\gm;\mathbb F_2)\to H^*_{et}(X\times \gm;\mathbb F_2)$ is an isomorphism for a finite type, smooth $X$ over $k$. The proof is much the same as that of Lemma 10.2 in \cite{milne}; the only observation we need here is that there is a quasi-isomorphism $f^*P^\bullet\simeq \mathbf Rp_*(q^*\mathbb F_2)$, where $P^\bullet$ is the complex $\mathbb F_2\stackrel 0\to\mathbb F_2 \to 0$, and $p$ and $q$ are projections from $X\times \gm$ on the first and second components, respectively. This is immediate from the Gysin sequence of the trivial $\gm$-bundle on $X$.

With this preparation, we observe that the cohomology classes in $H^*_{sing}(BGO(2n);\mathbb F_2)$ introduced in section 4.3 of Holla and Nitsure \cite{nitin1} make sense in the \'etale setup, and the proof of Theorem 4.3.5 of \cite{nitin1} carries over word by word. Therefore the description in terms of explicit generators and relations of the \'etale cohomology ring $H^*_{et}(BGO(2n);\mathbb F_2)$ of the classifying stack of $GO(2n)$ over any algebraically closed field is exactly the same as the reduced mod-2 singular cohomology ring of the topological space $BGO(2n)$, where $GO(2n)$ is considered as a Lie group over complex numbers. Also, the map induced in \'etale cohomology by the inclusion $BGO(2n)\subset BGL_{2n}$ has exactly the same description as in Proposition 3.2 of \cite{nitin2}.

\centerline{\bf Characteristic Classes for $S\!O(2n+1)$}
It is perhaps well known that $H^*(BS\!O(2n+1))\cong \mathbb F_2[w_2,\ldots,w_{2n+1}]$. Nonetheless, for the lack of a suitable reference in the cadre of algebraic stacks, we sketch a proof of this. Note that $O(2n+1)\cong \mu_2\times S\!O(2n+1)$, for $n\ge 1$. Hence, $BO(2n+1)\cong B \mu_2\times BS\!O(2n+1)$ for $n\ge 1$. We will show that in this specific case K\"unneth formula holds for the cohomology of algebraic stacks: $H^*(B\mu_2\times BS\!O(2n+1))\cong H^*(B\mu_2)\otimes H^*(BS\!O(2n+1))$. This will imply that $\mathbb F_2[w_1,\ldots,w_{2n+1}]\cong \mathbb F_2[w_1]\otimes H^*(BS\!O(2n+1))$, from which it can be seen that $H^*(BSO(2n+1))\cong\mathbb F_2[w_2,\ldots,w_{2n+1}]$. 

Claim: If $X$ is any finite type, smooth scheme over $k$, and if $B_r$ stands for the bundle quotient of $\mathbb A^{r+1}-\{0\}$ by $\mu_2\subset \gm$, then the natural map $H^*(B_r)\otimes H^*(X)\to H^*(B_r\times X)$ is an isomorphism. To prove it, we begin by observing the following: Let $E\to B$ be a principal $\gm$-bundle, and let $B'=E/\mu_2$, so that we have a natural projection $B'\to B$, which is again a principal $\gm$-bundle (because $\gm/\mu_2\cong\gm$). Then the class of the principal $\gm$-bundle $B'\to B$ in $H^1(B,\gm)$ is the square of that of $E\to B$, and therefore, its image under the connecting morphism coming from the Kummer sequence is zero in $H^2(B,\mu_2)$. This means that the Gysin sequence for $B'\to B$ splits into short exact sequences. 

Now look at the following commutative diagram, whose left column is a part of the Gysin sequence for $B_r\to\mathbb P^r$ tensored with $H^*(X)$, and the right column is the corresponding part of the Gysin sequence for $B_r\times X\to \mathbb P^r\times X$.

\centerline{\xymatrix{
0\ar[d]& 0\ar[d]\\
\oplus_{i+j=m}H^i(\mathbb P^r)\otimes H^j(X)\ar[r]^-\sim\ar[d] & H^m(\mathbb P^r\times X)\ar[d]\\
\oplus_{i+j=m}H^i(B_r)\otimes H^j(X)\ar[r]\ar[d] & H^m(B_r\times X)\ar[d]\\
\oplus_{i+j=m}H^{i-1}(\mathbb P^r)\otimes H^j(X)\ar[r]^-\sim\ar[d] & H^{m-1}(\mathbb P^r\times X)\ar[d]\\
0 & 0}}

\noindent The first and the third horizontal maps are isomorphisms, by Lemma 10.2 of \cite{milne}. Therefore the middle row is an isomorphism $H^*(B_r)\otimes H^*(X)\cong H^*(B_r\times X)$, as claimed.

The K\"unneth formula for $B\mu_2\times BS\!O(2n+1)$ can be deduced from this claim, as follows. As we have described in Section \ref{charstack}, we can use acyclic covers to approximate the cohomology of classifying stacks. In the formula $H^*(B_r)\otimes H^*(X)\cong H^*(B_r\times X)$, we substitute for $X$ an acyclic cover for $BS\!O(2n+1)$, and note that $B_r$ is an $r$-acyclic cover for $B\mu_2$. This completes the proof of the K\"unneth formula for $B\mu_2\times BS\!O(2n+1)$.
\vskip 2em
\noindent{\bf Acknowledgment:}
This result will form a part of the my PhD thesis. I thank my supervisor, Nitin Nitsure, for his guidance.

\noindent\emph{Address: School of Mathematics, Tata Institute of Fundamental Research, Homi Bhabha Road, Mumbai 400005, India.\\
e-mail: saurav@math.tifr.res.in}

\centerline{18-I-2012}
\end{document}